\documentclass[a4paper,fleqn]{article}
\newcommand{\floor}[1]{\lfloor{#1}\rfloor}
\newcommand{\ceil}[1]{\lceil{#1}\rceil}
\newcommand{\defeq}{\stackrel{\text{def}}{=}}
\newcommand{\F}[2]{{}_{#1}F_{#2}}

\usepackage{amsmath}

\title{The number of permutations with a given number of sequences}

\author{Marcus Kollar\\[1em]
   \small\parbox{0.7\textwidth}{\centering %
     \textit{Theoretical Physics III, Center for Electronic Correlations and Magnetism,\\
     Institute for Physics, University of Augsburg, 86135 Augsburg, Germany}\\
     \texttt{Marcus.Kollar@physik.uni-augsburg.de}}}

\date{\normalsize October 15, 2006}

\begin{document}
  \maketitle

  \begin{abstract}
    \vspace*{-2em}
    
    $P(n,s)$ denotes the number of permutations of~$1,2,\ldots n$ that
    have exactly~$s$ sequences.  Canfield and Wilf [math.CO/0609704]
    recently showed that $P(n,s)$ can be written as a sum of~$s$
    polynomials in~$n$.  We determine these polynomials explicitly and
    also obtain explicit expressions for $P(n,s)$ and its fixed-$s$
    generating function $u_s(x)$.
  \end{abstract}
  
  \section{Introduction}\label{sec:intro}
  
  A permutation of the numbers $1,2,\ldots n$ has at least two
  sequences, where a \textit{sequence} or \textit{run} is defined as
  the maximal interval of consecutive increasing or decreasing
  numbers.  Let $P(n,s)$ denote the number of permutations of $n$
  numbers with $s$ runs.  Its basic recurrence is
  \cite{Andre1884}%
  \begin{subequations}\label{eq:Peq}%
    \begin{align}%
      P(n,s)
      &=
      sP(n-1,s)
      +
      2P(n-1,s-1)
      +
      (n-s)P(n-1,s-2)
      \,,
      &
      (n&\geq 2)
      \\
      P(2,s)
      &=
      2\delta_{s,1}    
      \,,
      \,.
    \end{align}%
  \end{subequations}
  Recently Canfield and Wilf \cite{Canfield2006} showed that the
  fixed-$s$ generating function of $P(n,s)$ has the form
  \begin{align}
    u_s(x)
    &\defeq
    \sum_{n=2}^{\infty}P(n,s)x^n
    =
    \frac{\Phi_s(x)}{\Delta_s(x)}
    \,,
    &
    \Delta_s(x)
    &\defeq
    \prod_{i=0}^{s-1}(1-(s-i)x)^{\floor{i/2}+1}
    \,,\label{eq:ueq}
    &
    (s&\geq1)
  \end{align}
  where $\Phi_s(x)$ is a polynomial of degree $1+\ceil{s(s+2)/4}$,
  i.e., one larger than the degree of the polynomial $\Delta_s(x)$ in
  the denominator.  Furthermore they showed that $P(n,s)$ can be
  written as
  \begin{align}
    P(n,s)&=\sum_{i=0}^{s-1}\psi_i(n,s)(s-i)^n
    \,,\label{eq:PCW}
    &
    (n-1&\geq s\geq1)    
  \end{align}
  where each $\psi_i(n,s)$ is a polynomial in $n$ of degree at most
  $\floor{i/2}$.  These polynomials satisfy%
  \begin{subequations}\label{eq:psieq}%
    \begin{align}%
      \psi_i(n,s)
      &=
      \dfrac{
        s\psi_i(n-1,s)
        +
        2\psi_{i-1}(n-1,s-1)
        +
        (n-s)\psi_{i-2}(n-1,s-2)
      }{s-i}
      \,,
      &
      (s-1&\geq i\geq1)
      \\
      \psi_{-1}(n,s)&=0
      \,,
      \qquad
      \psi_{0}(n,s)=K(s)\defeq2^{-(s-2)}
      \,,
    \end{align}%
  \end{subequations}
  so that (\ref{eq:PCW}) satisfies (\ref{eq:Peq}).
  
  In Section~\ref{sec:polys} we determine the polynomials
  $\psi_i(n,s)$ from (\ref{eq:psieq}) and an explicit formula for
  $P(n,s)$ from (\ref{eq:PCW}). An auxiliary generating function is
  evaluated in Section~\ref{sec:Akz}. In Section~\ref{sec:usx} we
  obtain explicit expressions for the generating functions $u_s(x)$
  and the polynomials $\Phi_s(x)$.

  \section{The polynomials $\psi_i(n,s)$}\label{sec:polys}
  
  We put $\phi_i(n,t)\defeq\psi_i(n,s-i)$.
  {}From (\ref{eq:psieq}) we find%
  \begin{subequations}\label{eq:phieq}%
    \begin{align}%
      t\phi_i(n,t) &= (t+i)\phi_i(n-1,t) + 2\phi_{i-1}(n-1,t) +
      (n-t-i)\phi_{i-2}(n-1,t) \,,\label{eq:phirec} & (t&\geq 1,
      i\geq1)
      \\
      \phi_{-1}(n,t)&=0, \qquad \phi_{0}(n,t)=K(t)
      \,,\label{eq:phiini}
    \end{align}%
  \end{subequations}
  which has the advantage that $t$ is merely a parameter but not
  involved in the recurrence.  

  {}From (\ref{eq:phieq}) we find that the generating function
  $f(x,n,t)\defeq\sum_{i=0}^{\infty}\phi_i(n,t)x^i$ satisfies%
  \begin{subequations}\label{eq:feq}%
    \begin{align}%
      tf(x,n,t)
      &=
      \left(2x+(n-2)x^2+(1-x^2)t\right)f(x,n-1,t)
      +(1-x^2)x\frac{\partial}{\partial x}f(x,n-1,t)
      \,,\label{eq:frec}
      \\
      f(0,n,t)
      &=
      K(t)
      \,.\label{eq:fini}
    \end{align}%
  \end{subequations}
  We seek a solution of the form $f(x,n,t)=g(x,t)h(x,t)^n$. After
  inserting this into (\ref{eq:frec}) we first eliminate all terms
  proportional to $n$ by choosing $h(x,t)=(1-x^2)^{1/2}$.  This leads
  to a separable linear differential equation for $g(x,t)$, which is
  solved by $g(x,t)=(1-x)^{1/2}(1+x)^{-3/2}(1+(1-x^2)^{1/2})^{-t}$.
  Using (\ref{eq:fini}) we arrive at
  \begin{align}
    f(x,n,t)
    &=
    K(t)
    (1-x)^2
    (1-x^2)^{(n-3)/2}
    \left(\frac{2}{1+\sqrt{1-x^2}}\right)^{t}
    \,.\label{eq:fres}
  \end{align}
  In order to determine the coefficient of $x^i$ of (\ref{eq:fres}) we
  use \cite{Wilf1994}
  \begin{align}
    (1-x^2)^{(n-3)/2}
    &=
    \sum_{k=0}^{\infty}a_k(n)x^{2k}
    \,,
    &
    a_k(n)
    &\defeq
    \binom{(n-3)/2}{k}(-1)^k
    \,,
    \\
    \left(\frac{2}{1+\sqrt{1-x^2}}\right)^{t}
    &=
    \left(\frac{1-\sqrt{1-x^2}}{x^2/2}\right)^{t}
    =
    \sum_{m=0}^{\infty}b_m(t)x^{2m}
    \,,
    &
    b_m(t)
    &\defeq
    \frac{t(2m+t-1)!}{m!(m+t)!4^m}
    \,.
    \intertext{We find}
    f(x,n,t)
    &=
    K(t)
    (1-x)^2
    \sum_{j=0}^{\infty}p_j(n,t)x^{2j}
    \,,
    &
    p_j(n,t)
    &\defeq
    \sum_{k=0}^{j}a_k(n)b_{j-k}(t)
    \,.
  \end{align}
  With the convention $p_{-1}(n,t)=0$ we finally obtain
  \begin{align}
    \phi_{2j}(n,t)
    &=
    K(t)\left(p_j(n,t)+p_{j-1}(n,t)\right)
    \,,
    &
    \phi_{2j+1}(n,t)
    &=
    -2K(t)p_j(n,t)
    \,,\label{eq:phi}
    &
    (j&\geq0)
  \end{align}
  as well as 
  \begin{align}
    \psi_{i}(n,s)
    &=
    (-1)^i
    K(s-i)
    \left(p_{\big\lfloor\tfrac{i}{2}\big\rfloor}(n,s-i)+p_{\big\lfloor\tfrac{i-1}{2}\big\rfloor}(n,s-i)\right)
    \hspace*{-8ex}{}
    \\
    &=
    K(s-i)
    \sum_{j=0}^{\floor{i/2}}
    g_{i,j}
    p_j(n,s-i)
    \,,
    &
    g_{i,j}
    \defeq
    \begin{cases}
      \delta_{j,i/2}+\delta_{j,i/2-1}
      &
      i\text{ even}
      \\
      -2\delta_{j,(i-1)/2}
      &
      i\text{ odd}
    \end{cases}
    \,.\label{eq:psi}
  \end{align}
  Since $a_j(n)$ is a polynomial in $n$ of order $j$, so are
  $p_j(n,t)$, $\phi_{2j}(n,t)$, and $\phi_{2j+1}(n,t)$. Thus
  $\psi_{i}(n,s)$ is a polynomial in $n$ \textit{precisely} of degree
  $\floor{i/2}$; it was already shown in Ref.~\cite{Canfield2006} that
  this degree is \textit{at most} $\floor{i/2}$.  The first few
  $\psi_{i}(n,s)$ are given in Table~\ref{tab:psitab}.
  
  Inserting (\ref{eq:psi}) into (\ref{eq:PCW}) immediately yields an
  explicit formula for $P(n,s)$,
  \begin{align}
    P(n,s)
    &=
    \sum_{i=0}^{s-1}K(s-i)(s-i)^n
    \sum_{j=0}^{\floor{i/2}}g_{i,j}
    p_{j}(n,s-i)
    \,,\label{eq:Pexpl}
    &
    (n-1&\geq s\geq1)    
  \end{align}
  which will be used again in Section~\ref{sec:usx}.

  \begin{table}[p]
    \centering
    \begin{tabular*}{\textwidth}{r|l}
      \hline
      $i$&$\psi_i(n,s)$\\
      \hline
      $ -1$ & $0$
      \\
      $  0$ & $K(s)$
      \\
      $  1$ & $K(s-1)(-2)$
      \\
      $  2$ & $K(s-2)(-2n+s+8)/4$
      \\
      $  3$ & $K(s-3)(2n-s-3)/2$
      \\
      $  4$ & $K(s-4)(4n^2-4(s+8)n+s^2+15s+32)/32$
      \\
      $  5$ & $K(s-5)(-8n^2+8(s+3)n-2(s^2+5s+10))/32$
      \\
      $  6$ & $K(s-6)(-8n^3+12(s+8)n^2-2(3s^2+45s+98)n+s^3+21s^2+74s+144)/384$  
      \\
      $  7$ & $K(s-7)(16n^3-24(s+3)n^2+4(3s^2+15s+32)n-2(s^3+6s^2+23s+42)/384$
      \\
      $  8$ & $K(s-8)(16n^4-32(s+8)n^3+8(3s^2+45s+100)n^2-8(s^3+21s^2+76s+160)n$
      \\
            & $+s^4+26s^3+107s^2+442s+768)/6144$
      \\
      $  9$ & $K(s-9)(-16n^4+32(s+3)n^3-8(3s^2+15s+34)n^2+8(s+3)(s^2+3s+16)n$
      \\
            & $-(s^4+6s^3+35s^2+126s+216))/3072$
      \\
      $ 10$ & $K(s-10)(-32n^5+80(s+8)n^4-80(s^2+15s+34)n^3+40(21s^2+s^3+78s+176)n^2$
      \\
            & $-2(5s^4+130s^3+555s^2+2510s+4504)n$
      \\
            & $+s^5+30s^4+115s^3+870s^2+2824s+4800)/122880$
      \\
      \hline
    \end{tabular*}
    \caption{The first few polynomials $\psi_{i}(n,s)$.}
    \label{tab:psitab}
  \end{table}

  \section{The auxiliary generating functions $A_k(z)$}\label{sec:Akz}
  
  Before turning to $u_s(x)$ we evaluate the auxiliary generating
  functions $A_{k}(z)\defeq\sum_{n=2}^{\infty}a_k(n)z^n$.  Their
  generating function is
  \begin{align}
    A(y,z)
    &\defeq
    \sum_{k=0}^{\infty}
    A_k(z)y^k
    =
    \sum_{n=2}^{\infty}(1-y)^{(n-3)/2}z^n
    \nonumber\\
    &=
    \frac{1}{\sqrt{1-y}}\frac{z^2}{1-z\sqrt{1-y}}
    =
    \left(z+\frac{1}{\sqrt{1-y}}\right)
    \frac{z^2}{1-z^2(1-y)}
    \nonumber\\
    &=
    \sum_{k=0}^{\infty}
    \left(
      \frac{z^{2k+3}}{(1-z^2)^{k+1}}
      +
      \sum_{m=0}^{\infty}
      \binom{m-1/2}{k}z^{2m+2}
    \right)(-y)^k
    \,.
  \end{align}
  The inner sum evaluates to
  \begin{align}
    \sum_{m=0}^{\infty}
    \binom{m-1/2}{k}x^{m}    
    &=
    \binom{-1/2}{k}\F{2}{1}(\tfrac12,1;\tfrac12-k;x)
    =
    \binom{-3/2}{k}\frac{(-x)^k}{(1-x)^{k+1}}\F{2}{1}(-k,\tfrac12;\tfrac32;x^{-1})
    \,.
  \end{align}
  $A_k(z)$ is then obtained as
  \begin{align}
    A_k(z)
    &=
    \frac{z^2(-1)^k\widetilde{A}_k(z)}{(1-z^2)^{k+1}}
    \,,
    &
    \widetilde{A}_k(z)
    &\defeq
    z^{2k+1}
    +
    \binom{-3/2}{k}
    \sum_{m=0}^{k}
    \binom{k}{m}
    \frac{(-1)^{k+m}}{2m+1}z^{2k-2m}
    \,,\label{eq:Atilde}
  \end{align}
  where $\widetilde{A}_k(z)$ is a polynomial in $z$ of degree $2k+1$.
  We evaluate its $k$th derivative at $z=-1$ and find
  \begin{align}
    \left.\frac{d^k\widetilde{A}_k(z)}{dz^k}\right|_{z=-1}
    &=
    \binom{-3/2}{k}\frac{k!}{k+1}
    \left[
    -4^k
    +
    \sum_{m}
    \binom{k}{m}\binom{2k-2m}{k-2m}
    \frac{(-1)^m(k+1)}{2m+1}
    \right]
    =0
    \,,
  \end{align}
  since the sum which appears inside the square brackets equals $4^k$,
  as certified by the WZ \cite{Wilf1994} proof certificate $R(k,m)
  =4m(2m+1)(2k-2m+1)(k+1-2m)^{-1}(k+1)^{-1}$. Thus $z=-1$ is an at
  least $(k+1)$-fold zero of $\widetilde{A}_k(z)$, i.e., the
  polynomial $\widetilde{A}_k(z)$ contains a factor $(1+z)^{k+1}$.

  For the higher derivatives of $\widetilde{A}_k(z)$ at $z=-1$ we find
  \begin{align}
    \widetilde{a}_k(p)
    &\defeq
    \frac{(-1)^k}{(k+p+1)!}
    \left.\frac{d^{k+p+1}\widetilde{A}_k(z)}{dz^{k+p+1}}\right|_{z=-1}
    \\
    &=
    (-1)^{p}
    \left[
      \binom{2k+1}{k+p+1}
      -
      \binom{-3/2}{k}
      \sum_{m=0}^{\floor{(k-p-1)/2}}
      \binom{k}{m}\binom{2k-2m}{k+p+1}
      \frac{(-1)^{k+m}}{2m+1}
    \right]
    \,.
    &
    (0&\leq p\leq k)
  \end{align}
  The remaining sum (which vanishes for $p=k$) does not seem to allow
  further simplification.
  
  After inserting
  $\widetilde{A}_k(z)=(-1)^k\sum_{p=0}^{k}\widetilde{a}_k(p)(1+z)^{k+p+1}$
  into (\ref{eq:Atilde}), the final expression for $A_k(z)$ becomes
  \begin{align}
    A_k(z)
    &=
    \frac{\widetilde{\Phi}_k(z)}{(1-z)^{k+1}}
    \,,
    &
    \widetilde{\Phi}_k(z)
    &=
    z^2\sum_{p=0}^{k}\widetilde{a}_k(p)(1+z)^{p}
    \,.
    &
    (k&\geq0)
  \end{align}
  Note that the polynomial $\widetilde{\Phi}_k(z)$ in the numerator of
  $A_k(z)$ has degree $k+2$, one larger than the degree of the
  polynomial in the denominator.

  \section{The generating functions $u_s(x)$ and polynomials $\Phi_s(x)$}\label{sec:usx}
  
  For the generating functions $u_s(x)$ we employ (\ref{eq:Pexpl}) and
  $\sum_{n=2}^{\infty}p_j(n,t)z^n$ $=$
  $\sum_{k=0}^{j}A_k(z)b_{j-k}(t)$, whence
  \begin{align}
    u_s(x)
    &=
    \sum_{i=0}^{s-1}K(s-i)
    \sum_{j=0}^{\floor{i/2}}
    g_{i,j}
    \sum_{k=0}^{j}A_k((s-i)x)b_{j-k}(s-i)
    \,,
    ~~~~~~~~~
    (s\geq1)
    \nonumber\\
    &=
    \sum_{i=0}^{s-1}
    \sum_{k=0}^{\floor{i/2}}
    \frac{B_{i,k}(x,s-i)}{(1-(s-i)x)^{k+1}}
    \,,
    ~~~~~~~~~    
    B_{i,k}(x,t)
    \defeq
    K(t)
    \widetilde{\Phi}_k(tx)
    \sum_{j=k}^{\floor{i/2}}
    g_{i,j}
    b_{j-k}(t)
    \,.
  \end{align}
  Clearing denominators yields $u_s(x)=\Phi_s(x)/\Delta_s(x)$ as in
  (\ref{eq:ueq}), with
  \begin{align*}
    \Phi_s(x)
    =
    \sum_{i=0}^{s-1}
    \left(
      \prod_{\substack{m=0\\m\neq i}}^{s-1}(1-(s-i)x)^{\floor{m/2}+1}
    \right)
    \sum_{k=0}^{\floor{i/2}}
    B_{i,k}(x,s-i)
    (1-(s-i)x)^{\floor{i/2}-k}
    \,.
  \end{align*}
  The polynomial $\Phi_s(x)$ indeed has degree $1+\ceil{s(s+2)/4}$, as
  proven already in \cite{Canfield2006}. The first few $\Phi_{s}(x)$
  are given in Table~\ref{tab:Phitab}.
  
  \begin{table}[p]
    \centering
    \begin{tabular*}{\textwidth}{r|l}
      \hline
      $s$&$\Phi_s(x)$\\
      \hline
      $  1$ & $2x^2$
      \\
      $  2$ & $4x^3$
      \\
      $  3$ & $2x^4(5-6x)$
      \\
      $  4$ & $4x^5(24x^2-29x+8)$
      \\
      $  5$ & $2x^6(720x^4-1704x^3+1436x^2-501x+61)$
      \\
      $  6$ & $4x^7(17280x^6-51336x^5+61188x^4-37256x^3+12209x^2-2041x+136)$
      \\
      $  7$ & $2x^8(-3628800x^9+15729120x^8-29341872x^7+30810864x^6-20028656x^5+8353808x^4$\\
      &$-2236439x^3+370871x^2-34601x+1385)$
      \\
      $  8$ & $4x^9(696729600x^{12}-3555239040x^{11}+8107966944x^{10}-10906662240x^9+9627417336x^8$\\
      &$-5872225480x^7+2537780728x^6-783164808x^5+171355239x^4-25936503x^3+2579241x^2$\\
      &$-151385x+3968)$
      \\
      $  9$ & $2x^{10}(1316818944000x^{16}-8712886694400x^{15}+26410986334080x^{14}-48618945021312x^{13}$\\
      &$+60779114417952x^{12}-54684478479456x^{11}+36624658707312x^{10}-18628018251952x^9$\\
      &$+7273896122392x^8-2188789058612x^7+506111568077x^6-89028957282x^5+11685816855x^4$\\
      &$-1107016832x^3+71414171x^2-2804314x+50521)$
      \\
      $ 10$ & $4x^{11}(2528292372480000x^{20}-18993012627456000x^{19}+66507291476582400x^{18}$\\
      &$-144199874533248000x^{17}+216971940209451264x^{16}-240735551604776064x^{15}$\\
      &$+204330019791468672x^{14}-135856983272339904x^{13}+71875337579512880x^{12}$\\
      &$-30562090468050280x^{11}+10504633067351272x^{10}-2924633644527940x^9+658629666786430x^8$\\
      &$-119364099863329x^7+17244871619376x^6-1956223222079x^5+170214919190x^4$\\
      &$-10952481287x^3+490431140x^2-13630637x+176896)$
      \\
      \hline
    \end{tabular*}
    \caption{The first few polynomials $\Phi_s(x)$.}
    \label{tab:Phitab}
  \end{table}

  \section{Conclusion}\label{sec:concl}
  
  We derived explicit expressions for $P(n,s)$ and $u_s(x)$ using the
  formulation of Ref.~\cite{Canfield2006} in terms of polynomials
  $\psi_i(n,s)$ and $\Phi_s(x)$.  Previously, Carlitz derived explicit
  formulas for $P(n,s)$ based on a two-variable generating function
  \cite{Carlitz1978,Carlitz1980,Carlitz1981}.  However, Canfield and
  Wilf \cite{Canfield2006} note that the final formulas in
  Ref.~\cite{Carlitz1980} are not entirely correct.  We checked our
  results by computer algebra for $i,s,n\leq20$.
  
  \vspace*{\fill}

\end{document}